\documentclass[prl,aps]{revtex4}
\usepackage{graphicx}

\begin{document}

\title[A lattice model related to the nonlinear Schroedinger equation]{A lattice model related to the nonlinear Schroedinger equation}
 {\rm Doklady  Akademii  Nauk vol  259, page 76, l981 }
\bigskip

\author{A.\ G.\ Izergin and V.\ E.\ Korepin}
\affiliation{Leningrad branch  of  Steklov  Mathematical  Institute. of Academy of Sciences  of the USSR}
\date{Feb 26, 1981}

\begin{abstract}
This is a historical note.  In 1981 we constructed a discrete version of quantum  nonlinear Schroedinger equation.
This led to our discovery of quantum determinant: it appeared in construction of anti-pod (11). Later these became important in
quantum groups: it describes center of Yang-Baxter algebra.  Our paper was published in Doklady  Akademii  Nauk vol  259, page 76 (July l981) in Russian language.  
\end{abstract}

%\pacs{Pacs No: 03.67.-a, 03.65.Ud, 32.80.Lg}

\maketitle

\begin{enumerate}

 \item
 Inverse scattering method is used  (in classical \cite{ggkm, zafa} and quantum cases \cite{fad}) to solve evolutionary equations of completely integrable dynamical systems.  In quantum case we shall abbreviate the method to qism. It is based on Lax representation 
 \begin{equation} \label{lax}
\partial_{t}L_{n}(\lambda)= M_{n+1}(\lambda)L_{n}(\lambda) - L_{n}(\lambda) M_{n}(\lambda)
 \end{equation}
 Entries  of the matrices $L_{n}(\lambda) $ and $M_{n}(\lambda) $
are expressed in terms of dynamical variables of the lattice system (they also depend on spectral parameter $\lambda$).
 The monodromy matrix  $T^{n,m}(\lambda)=L_{n}(\lambda) \ldots L_{m}(\lambda)$  ($n\ge m$) and transfer matrix  $\tau(\lambda) =\mbox{tr}  T^{N,1}(\lambda) $ are important for qism.
Recently an effective method of construction of  action-angle
variables was discovered, it is based on $R$-matrix. In classical case \cite{skl1} it provides Poisson brackets of matrix elements of
the monodromy matrix, while in quantum case \cite{fad,bax} it gives commutation relations: 
\begin{eqnarray}  
& \left[ T_{c}(\lambda) \otimes_{\hskip -.18cm,} {\hskip .1cm}T_{c}(\mu)  \right] =  \left[ T_{c}(\lambda) \bigotimes T_{c}(\mu) , R_{c} (\lambda , \mu ) \right] \\
&R_{q}(\lambda , \mu ) \left( T_{q}(\lambda) \bigotimes T_{q}(\mu) \right) = \left( I \bigotimes T_{q}(\mu) \right) \left( T_{q}(\lambda) \bigotimes I \right) R_{q}(\lambda , \mu )
\end{eqnarray} 
These equations lead to commutativity of transfer matrices (we use subindex c or q to distinguish classical from quantum).  
Meaning that $\partial_{\mu} \ln \tau (\mu)$ is a generating
functional of  Hamiltonians for completely integrable
systems.

In classical case E.K. Sklyanin proved \cite{skl2} that corresponding equations of motion can be represented in the Lax form (\ref{lax}).  Here we prove that this is true in quantum case as well. Generating functional of the operators $M_{n}(\lambda)$ is a matrix
$m_{n}(\lambda , \mu)$:
\begin{eqnarray}
& m_{n}(\lambda , \mu) = i \tau^{-1}(\mu)  \partial_{\mu}  \tau (\mu) -  i q^{-1}_{n}(\lambda , \mu) \partial_{\mu} q_{n}(\lambda , \mu) \nonumber \\
& q_{n}(\lambda , \mu) =  \mbox{tr}_{2} \left(  I \bigotimes T^{N,n} (\mu) \right) R_{q}^{-1}(\lambda , \mu)    \left(I \bigotimes T^{n-1,1} (\mu)     \right) \nonumber  \\
& i\left[ \partial_{\mu} \ln \tau (\mu), L_{n}(\lambda) \right] = m_{n+1} (\lambda , \mu) L_{n}(\lambda) - L_{n}(\lambda) m_{n} (\lambda , \mu) \nonumber
\end{eqnarray} 
Here $\mbox{tr}_{2}$ denotes trace in the second linear space of the tensor product.  The proof follows from
$q_{n+1}(\lambda , \mu)  L_{n}(\lambda) =  L_{n}(\lambda) q_{n}(\lambda , \mu)  $ .  So we proved that even in the quantum case 
it is sufficient to have $L_{n}(\lambda)$ operator and $R$ matrix in order to apply QISM.
\item
Let consider nonlinear Schroedinger equation (nS). In continuous case it has a Hamiltonian
\begin{eqnarray} \label{contham}
& H=  \int dx  \left( \partial_{x} \psi^{\dagger}  \partial_{x} \psi +\kappa \psi^{\dagger} \psi^{\dagger} \psi \psi \right), \nonumber \\
& \{\psi_{c} (x),  \psi_{c}^{\dagger} (y)  \} = i \delta (x-y), \qquad \left[ \psi_{q} (x), \psi_{q}^{\dagger} (y)  = \delta (x-y)\right]
\end{eqnarray} 
It is integrable both in classical \cite{zasha} and quantum \cite{bepo, fad} cases. Corresponding $R$- matrix can be called quasi-classical 
\begin{equation}
 R_{q}=I\bigotimes I -i R_{c}, \qquad   R_{c}=\frac{\kappa \Pi}{\lambda- \mu}
 \end{equation}
 Here $I$ is identical matrix 2X2 and $\Pi$ is parmutation matrix.
 
 Lattice generalization of nS has long attracted attention of the experts \cite{abl,kul}. We propose a new version of lattice nS both in classical and quantum cases. It is distinguishing feature is that $R$-matrix is the same as in the continuous case and basic variables $\chi$ are canonical Bose fields. We start by suggesting the following $L_{n}$ operator:
 \begin{eqnarray} \label{spin}
 & L_{n}(\lambda) = -i \lambda \Delta \sigma_{3}/2 + S^{3}_{n} I +S^{+}_{n}\sigma_{+} + S^{-}_{n}\sigma_{-} ,       \nonumber \\
 & S^{3}=1+\frac{\kappa}{2} \chi^{\dagger}_{n} \chi_{n} , \qquad S^{+}_{n}= -i\sqrt{\kappa} \chi^{\dagger}_{n} \rho_{n}^{+},
 \qquad  S^{-}_{n}= i\sqrt{\kappa} \rho_{n}^{-} \chi_{n}\nonumber \\
 & \rho_{n}^{\pm}= \rho_{n}^{\pm} (\chi_{n} \chi_{n}^{\dagger}  ) , \qquad  \rho_{n}^{+}  \rho_{n}^{-}= 1+\frac{\kappa}{4} \chi^{\dagger}_{n} \chi_{n} , \qquad 2\sigma_{\pm}= \sigma_{1}\pm i\sigma_{2} \nonumber \\
  & \{ \chi^{c}_{m}, {\chi^{c}_{n}}^{\dagger} \} =i\Delta \delta_{m,n}, \qquad  \left[ \chi^{q}_{m}, {\chi^{q}_{n}}^{\dagger} \right] =\Delta \delta_{m,n}
  \end{eqnarray}
Here $\sigma$ are Pauli matrices. We consider repulsive case $\kappa >0$ and put $\rho^{+}_{n}=\rho^{-}_{n}= \rho_{n} $.
\item 
Here we shall discuss  lattice model (\ref{spin}) in classical case. Simple calculations lead to
\begin{eqnarray} \label{cldet}
&T(\lambda)\sigma_{2}T^{t} (\lambda) \sigma_{2} = d_{c}^{n-m+1} (\lambda) I \nonumber \\
& d_{c} (\lambda)= \mbox{det}  L_{n}(\lambda) = 1+ \lambda^{2} \Delta^{2}/4 
\end{eqnarray}
This shows that at  $\lambda= \nu=-2i/\Delta$ the $L_{n} (\lambda)$ operator turns into one dimensional projector.
This makes it possible to calculate explicitly logarithmic derivatives of $\tau (\lambda)$ at this point, which can be represented as a sum of local densities:
\begin{eqnarray} \label{clogder}
&\partial_{\mu}^{n}\ln \tau (\lambda)|_{\lambda=\nu}= \sum_{k=1}^{N}  h_{k,n} , \\
&h_{k,n}=D^{n}\ln \mbox{tr} L_{k+n}(\nu) L_{k+n-1}(\lambda_{k+n-1})  \ldots L_{k}(\lambda_{k}) L_{k-1}(\nu)|_{\lambda_{.}=\nu}  \nonumber
\end{eqnarray}
Here $D^{n}$ is a differential operator. For small $n$ it is:
\begin{eqnarray} \label{diffoper}
& D^{1}=\partial_{k}=\frac{d}{d \lambda_{k}}, \qquad  D^{2}= 2\partial_{k+1}\partial_{k} + \partial^{2}_{k}  \\
&  D^{3}= 6\partial_{k+2}  \partial_{k+1} \partial_{k} + 6 \partial_{k+2}^{2}\partial_{k+1} + 6 \partial_{k+2} \partial_{k+1}^{2}-6\partial_{k+2}^{2}\partial_{k} -6 \partial_{k+2}\partial_{k}^{2}  + \partial_{k}^{3}\nonumber
\end{eqnarray}
We use this notations to define lattice classical Hamiltonian of nS
\begin{eqnarray} \label{latham}
& H_{c}=D_{c}(\lambda) \ln \left[ \left( 1+\lambda/\nu  \right)^{-N}\tau (\lambda)   \right]  +\mbox{complex conjugate}\nonumber \\
& D_{c}(\lambda) = \frac{i}{12\kappa} \left( \frac{d}{d \lambda^{-1}} \right)^{3}
\end{eqnarray}
The explicit expression shows that this Hamiltonian describes interaction of  five nearest neighbors  on the lattice. 
In the continuous limit [$\chi_{n}= \psi_{n} \Delta$; $\psi_{n+1}-\psi_{n}=O(\Delta)$; $\Delta \rightarrow 0$, $N\rightarrow \infty$ but $N\Delta =$ const]
it goes to the correct  Hamiltonian of the continuous model (\ref{contham}) and $L_{n}$ operator (\ref{spin}) turns into correct 
 continuous $L_{n}$ operator, see \cite{zasha}.
\item
Here we construct quantum lattice nS model.  Quantum analog of  (\ref{cldet}) is given by
\begin{eqnarray} \label{qdet}
&T(\lambda)\sigma_{2}T^{t} (\lambda +i \kappa) \sigma_{2} = d_{q}^{n-m+1} (\lambda)I \nonumber \\
& d_{q} (\lambda) = \Delta^{2}(\lambda - \nu)(\lambda - \nu +i\kappa  )/4 
\end{eqnarray}
This defines {\bf {quantum determinant}}:
$$\mbox{det}_{q}T(\lambda)=T_{11}(\lambda)T_{22}(\lambda+i\kappa) - T_{12}(\lambda) T_{21}(\lambda +i\kappa)=d_{q}^{n-m+1} (\lambda)$$
To define Hamiltonian of the model let us add quantum correction like in \cite{fad}:
\begin{equation}
H_{q}=\left(  D_{c}(\lambda) +\frac{i\kappa}{6} \frac{d}{d \lambda^{-1}} \right) \ln \left[ \left( 1+\lambda/\nu  \right)^{-N}\tau (\lambda)   \right]  +\mbox{hermitian conjugate}
\end{equation}
The model can be solved by qism \cite{fad}.  The pseudo-vacuum $\Omega$ is annihilated by lattice Bose fields
$\chi_{n} \Omega =0$. The eigenvectors are given by algebraic Bethe ansatz:
$$\Psi (\lambda_{1} \ldots \lambda_{n})= B  (\lambda_{1} ) \ldots B  (\lambda_{n} ) \Omega , \qquad B  (\lambda )= T_{12}(\lambda)$$
These $\lambda_{j}$ satisfy a system of Bethe equations:
\begin{equation}
\left(  \frac{1-i\lambda_{j}\Delta /2}{1+i\lambda_{j}\Delta /2} \right)^{N}=\prod_{k\neq j} \frac{\lambda_{j}-\lambda_{k}-i\kappa}{\lambda_{j}-\lambda_{k}+i\kappa} \nonumber
\end{equation}
Corresponding eigenvalue of $\tau(\lambda)$ is
\begin{equation}
\left( 1-\frac{i\lambda \Delta}{2}  \right)^{N} \prod_{k=1}^{n}\frac{\lambda -\lambda_{k}+i\kappa}{\lambda -\lambda_{k}} +
\left( 1+\frac{i\lambda \Delta}{2}  \right)^{N} \prod_{k=1}^{n}\frac{\lambda_{k} -\lambda+i\kappa}{\lambda_{k} -\lambda} \nonumber
\end{equation}
From here we obtain energy levels [eigenvalues of the Hamiltonian]
\begin{eqnarray} \label{qenergy}
& H_{q} \Psi =\left( \sum_{k=1}^{n}E(\lambda_{k})  \right)  \Psi , \qquad E(\mu)= f(\mu)+{ \overline{f({\overline\mu)}} }\nonumber \\
& f(\mu)= \left( D_{c} + \frac{i\kappa}{6} \frac{d}{d \lambda^{-1}} \right) \ln  \left( \frac{\mu -\lambda +i\kappa}{\mu -\lambda}  \right)|_{\lambda=\nu} \nonumber
\end{eqnarray}
The Hamiltonian has correct   continuous limit (\ref{contham}) and $E(\mu)\rightarrow \mu^{2}$.
\item
Quantum nS model constructed above can be considered as a generalization of XXX model with negative spin $-2/\kappa \Delta$.
We can rewrite the $L_{n}$ operator (\ref{spin}) in the way similar to XXX:
$$L^{X}_{n}=-\sigma_{3}L_{n}=i\lambda +t^{k}_{n}\otimes \sigma_{k} $$
Here $t^{k}_{n}$ are simple linear combinations of $S^{k}_{n}$ from (\ref{spin}). They form 
an infinite dimensional  representation of $SU(2)$ algebra, see  \cite{hp}.
\end{enumerate}


\begin{thebibliography}{99}

\bibitem{ggkm}  C.S. Gardner, J.M. Green, M.D. Kruskal and R. M. Miura,  Phys. Rev. Lett. vol 19, page 1095, 1967
\bibitem{zafa}  V.E. Zakharov and L.D. Faddeev, Functional  Analyse and Applications vol 5, page 280, 1971
\bibitem{fad} L.D. Faddeev, Preprint R-2-79, LOMI, Acad, Sc, USSR, 1979
\bibitem{skl1} E.K. Sklyanin, Preprint R-3-79, LOMI, Acad, Sc, USSR, 1979
\bibitem{bax} R.J. Baxter. Ann. Phys. vol 70, page 1, 1973 
\bibitem{ik} A.G.Izegin and V.E. Korepin Preprint E-3-80, LOMI, Acad, Sc, USSR, 1980
\bibitem{skl2} E.K. Sklyanin, Zap. Nauch. Semin. LOMI, 1980
\bibitem{zasha} V.E. Zakharov and  A.B. Shabat. Sov. Phys. JETP, vol 34, page 62, 1972
\bibitem{bepo} F.A. Berezin, G.P. Pokhil and V. M. Finkel'berg, Vestnik Mosk. Univer. Ser 1 no 1 page 21, 1964 
\bibitem{sklfa} E.K. Sklyanin and L.D. Faddeev, Doklady  Akademii  Nauk vol   243, 1978 (Sov. Phys. Dokl vol 23, page 902, 1978)
\bibitem{skl3} E.K. Sklyanin, Doklady  Akademii  Nauk vol 244, page 1337 (Sov. Phys. Dokl vol 24, page 107, 1979)
\bibitem{abl}  M.Ablowitz, Stud. Appl. Math. vol 58, page 17, 1978
\bibitem{kul}  P.P. Kulish, Lett Math. Phys vol 5 page 111, 1981
\bibitem{hp} T.Holstein and H. Primakoff, Phys. Rev vol 58, page 1098, 1940
\end{thebibliography}
\end{document}